\documentclass{amsart}
\usepackage{amssymb}
\usepackage{amsthm}
\usepackage{amscd}
\usepackage{graphics}

\newtheorem{introprop}{Proposition}
\newtheorem*{introprop1A}{Proposition 1A}
\newtheorem*{introprop2A}{Proposition 2A}
\newtheorem*{introprop3A}{Proposition 3A}
\newtheorem*{introprop4A}{Proposition 4A}
\newtheorem*{introprop2B}{Proposition 2B}

\newtheorem{lemma}{Lemma}
\newtheorem{smlemma}{Lemma}[section]

\newtheorem{prop}{Proposition}[section]

\newtheorem*{ques}{Question}
\numberwithin{equation}{section}
\title{Isoperimetric inequalities and rational homotopy invariants}
\author{Larry Guth}
\address{Department of Mathematics, Stanford, Stanford CA, 94305 USA}
\email{lguth@math.stanford.edu}

\begin{document}
\begin{abstract} 
We estimate the linear isoperimetric constants of an
n-dimensional ellipse.  Using these estimates and a technique
of Gromov, we estimate the Hopf and linking invariants of
Lipschitz maps from ellipses to round spheres.  Using
these estimates, we give a lower bound for the k-dilation
of degree non-zero maps between ellipses.
\end{abstract}

\maketitle

We will discuss the following question relating geometry
and topology.

\begin{ques} Suppose that $F$ is a map from $(S^3, g)$ to $(S^2, h)$.  
If $F$ has Lipschitz constant $L$, how big can the 
Hopf invariant of $F$ be?  In particular, what is the asymptotic behavior
as $L \rightarrow \infty$?
\end{ques}

In this paper, we survey what is known about this question.  We begin
with some general results due to Gromov and then work out some examples,
especially the example of 3-dimensional ellipses.  After that, we consider some
generalizations of this question.  The paper contains a mix of exposition
and new work.

Gromov looked at this question in \cite{Gr1} and \cite{Gr3} and proved
two basic results.  The first result deals with the case that
$g$ and $h$ are the standard round metrics.

\begin{introprop} (Gromov, \cite{Gr3}, page 358) If $F$ is a map from $(S^3, stan)$ to $(S^2, stan)$
with Lipschitz constant $L$, then the Hopf invariant of $F$ has norm
at most $C L^4$.  On the other hand, for any $L > 2$, there exist maps with Lipschitz
constant $L$ and Hopf invariant at least $c L^4$.
\end{introprop}

Gromov generalized
this result to deal with other metrics $(S^3, g)$, giving an estimate
involving isoperimetric information about the metric $g$.  To state
Gromov's estimate, we need a little vocabulary.  If $z$ is a closed integral 1-cycle
in $(S^3, g)$, the filling volume of $z$ is the smallest volume of any
integral 2-chain $y$ with $\partial y = z$.  The isoperimetric constant
$Iso(g)$ is the supremal ratio $FillVol(z) / Vol(z)$ as $z$ varies over all
closed 1-cycles in $(S^3, g)$.

\begin{introprop} (Gromov, \cite{Gr1}, page 96) Let $F: (S^3, g)
\rightarrow (S^2, h)$ be a map with Lipschitz constant $L$.  Then
the Hopf invariant of $F$ obeys the following inequality.

$$|\textrm{Hopf }(F)| \le Iso(g) Vol(g) Area(h)^{-2} L^4. \eqno{(*)}$$

\end{introprop}

Given a metric $(S^3, g)$, it is not trivial to estimate $Iso(g)$.  The main
new result in this paper is an estimate for $Iso(g)$ when $g$ is the metric
of a 3-dimensional ellipse.

We recall that an n-dimensional ellipse with principal axes
$E_0 \le ... \le E_n$ is the set $\{ x \in \mathbb{R}^{n+1} | \sum_{i=0}^n
(x_i / E_i)^2 = 1 \}$.  An n-dimensional ellipse is one of the simplest
examples of a Riemannian metric on $S^n$.

\begin{introprop} Suppose that $E$ is a 3-dimensional ellipse with
principal axes $E_0 \le ... \le E_3$.  Then, up to a factor of 100, $Iso(E) \sim
E_2$.
\end{introprop}

If $g$ is the metric of an ellipse $E$ and $h$ is the standard metric, then the upper bound $(*)$
tells us that the Hopf invariant of $F$ has norm at most $C L^4 Vol(E) E_2 \sim
C L^4 E_1 E_2^2 E_3$.  This upper bound is essentially sharp.

\begin{introprop} Suppose that $g$ is an ellipsoidal metric on $S^3$ with principal
axes $E_0 \le ... \le E_3$.  If $L > C E_1^{-1}$, then there is a map $F$ from 
$(S^3, g)$ to $(S^2, stan)$ with Hopf invariant at least $c L^4 E_1 E_2^2 E_3$.
\end{introprop}

These matching upper and lower bounds give us an accurate idea of how the dimensions of the ellipse
$E$ relate to the largest Hopf invariant of a Lipschitz map from $E$.

In general the upper bound $(*)$ may be far from sharp.  We contrive an example where
$(*)$ is far from sharp by gluing together two differently shaped ellipses.

Next we turn to generalizations of this question in higher dimensions.  The Hopf
invariant is defined for maps from $S^{4n-1}$ to $S^{2n}$ for each integer $n \ge 1$.
The propositions above have analogues in every dimension.  But I have decided
instead to focus on the linking invariant.  The linking invariant is a homotopy invariant
for maps from $S^n$ to a wedge of spheres $S^{k_1} \vee S^{k_2}$.  It is well-defined
whenever the dimensions satisfy $n+1 = k_1 + k_2$ and each $k_i$ is at least 2.  The
linking invariant is defined in the following way.  Let $q_1$ be a generic point in $S^{k_1}$
and let $q_2$ be a generic point in $S^{k_2}$.  Consider the inverse images $F^{-1}(q_1)$
and $F^{-1}(q_2)$ in $S^n$.  These are each oriented submanifolds of $S^n$, and the condition
on the dimension guarantees that the two submanifolds have a well-defined linking number.
This number is the linking invariant of $F$.  The linking invariant is closely related to
the Hopf invariant, but there are more linking invariants to play with, because $n$ can take
any value and for each $n$ we have several choices of $k_1$ and $k_2$.

All of our previous results generalize to linking invariants.  (Gromov's results generalize
to all rational homotopy invariants, as discussed in Chap. 7 of \cite{Gr3}.)
For the rest of this introduction, we fix $n$, $k_1$, and $k_2$, with the convention that
$k_1 + k_2 = n + 1$ and $2 \le k_1 \le k_2$.

\begin{introprop1A} (Gromov) Equip $S^n$ and $S^{k_1} \vee S^{k_2}$ with the standard metrics.  Suppose
that $F: S^n \rightarrow S^{k_1} \vee S^{k_2}$ has Lipschitz constant $L$.  Then the linking
invariant of $F$ has norm at most $C(n) L^{n+1}$.  If $L$ is larger than $C(n)$, then there are
examples with linking invariant at least $c(n) L^{n+1}$.
\end{introprop1A}

To state the higher-dimensional version of Proposition 2, we need a little more vocabulary.
Suppose that $z$ is an exact integral k-cycle in a Riemannian manifold $(M^n, g)$.  The
filling volume of $z$ is the smallest volume of any integral (k+1)-chain $y$ with $\partial y = z$.
We define the isoperimetric constant $Iso_k(M^n, g)$ to be the supremal ratio $FillVol(z) / Vol(z)$
as $z$ varies over all exact k-cycles in $M$.

\begin{introprop2A} (Gromov) Equip $S^n$ with a metric $g$, and equip $S^{k_1} \vee S^{k_2}$ with
a metric $h_1$ on $S^{k_1}$ and $h_2$ on $S^{k_2}$.  Suppose that $F: S^n \rightarrow S^{k_1} \vee
S^{k_2}$ has Lipschitz constant $L$.  Then the linking invariant of $F$ is bounded as follows.

$$|L(F)| \le Iso_{n-k_1}(S^n, g) Vol(g) Area(h_1)^{-1} Area(h_2)^{-1} L^{n+1}.$$

\end{introprop2A}

In order to apply this upper bound, we have to estimate the isoperimetric constant of $(S^n, g)$.
This is not trivial in general.  In this paper we carry it out for n-dimensional ellipses.

\begin{introprop3A} Let $E$ denote the n-dimensional ellipse with principal
axes $E_0 \le ... \le E_n$.  Up to a constant factor $C(n)$, the 
isoperimetric constant $Iso_k(E)$ is given by the following
formula.

$$Iso_k(E) \sim E_{k+1} + E_{n-k}.$$

\end{introprop3A}

Proposition 3A is the main new result in this paper.  To get a sense for it,
let's consider some k-cycles $z$ in $E$.  First we consider the smallest
k-dimensional equator of $E$.  This equator is the intersection of $E$ with
the (k+1)-plane spanned by the $x_0, ..., x_k$ coordinates.  This equator
is a k-dimensional ellipse with principal axes $E_0 \le .. \le E_k$ and
volume roughly $E_1 ... E_k$.  Suppose now that we intersect $E$ with the 
(k+2)-plane spanned by the first k+2 coordinates.  The intersection is a (k+1)
dimensional ellipse with volume roughly $E_1 ... E_{k+1}$, and $z$ is
a hypersurface in the intersection, dividing it into two equal pieces.
Therefore, the filling volume of $z$ is at most $\sim E_1 ... E_{k+1}$.  But
it's not too hard to convince oneself that this filling is near-optimal,
so that the filling volume of $z$ really is $\sim E_1 ... E_{k+1}$.  In
this case $FillVol(z) / Vol(z)$ is roughly $E_{k+1}$.

Second, we consider the largest k-dimensional equator of $E$.  This equator
is the intersection of $E$ with the (k+1)-plane spanned by the $x_{n-k},
..., x_n$ coordinates.  It is an ellipse with principal axes $E_{n-k},
..., E_n$ and volume roughly $E_{n-k+1} ... E_n$.  By the same argument
as above, it bounds a (k+1)-dimensional hemisphere with volume $\sim
E_{n-k} ... E_n$, and it's not hard to see that its filling volume
is $\sim E_{n-k} ... E_n$.  In this case, the ratio $FillVol(z) / Vol(z)$
is roughly $E_{n-k}$.

This argument shows that $Iso_k(E)$ is at least $c(n) [E_{k+1} + E_{n-k}]$.
The proposition claims that $Iso_k(E)$ is not much bigger than this.  In other
words, the proposition tells us that an arbitrary cycle $z$ in $E$ is not
substantially harder to fill than one of the equators of $E$.

Combining Proposition 2A and Proposition 3A, 
we get an upper bound for the linking invariant
of a Lipschitz map from $E$, and the upper bound turns out to be
essentially sharp.

\begin{introprop4A} Let $E$ be an n-dimensional ellipse with principal
axes $E_0 \le ... \le E_n$.  Equip
$S^{k_1} \vee S^{k_2}$ with the standard metrics.  Suppose that $F$ is a map
from $E$ to $S^{k_1} \vee S^{k_2}$ with Lipschitz constant $L$.  Then the linking
invariant of $F$ has norm at most $C(n) L^{n+1} E_{n-k_1+1} E_1 ... E_n$.
On the other hand, if $L > C(n) E_1^{-1}$, then there exists a map with Lipschitz
constant $L$ and linking
invariant at least $c(n) L^{n+1} E_{n-k_1+1} E_1 ... E_n$.
\end{introprop4A}

I became interested in this subject because of my work on k-dilation.  The
k-dilation is a generalization of the Lipschitz constant which measures
how much a mapping stretches k-dimensional volumes.  If $F$ is a smooth
map between Riemannian manifolds, then $Dil_k(F)$ can be defined as the
supremum of $| \Lambda^k dF |$.  It can also be described more geometrically
as follows.  We say that $F$ has k-dilation at most $\Lambda$ if, for every
k-dimensional surface $\Sigma^k$ in the domain, the volume of $F(\Sigma^k)$
is at most $\Lambda$ times the volume of $\Sigma^k$.  If $F$ has Lipschitz
constant $L$, then it follows easily that it also has k-dilation at most
$L^k$ for every $k$.  On the other hand, a map with 2-dilation 1 may have
arbitrarily large Lipschitz constant.  A map with 2-dilation 1 may stretch
lengths by an arbitrarily large amount in one direction as long as it compensates
by contracting in other directions.

Gromov was one of the first mathematicians to look at k-dilation.  In particular,
he noticed that the proofs of the propositions above don't need the hypothesis
$Lip(F) \le L$, but only the hypothesis $Dil_k(F) \le L^k$ for an appropriate
value of $k > 1$.  For example, Proposition 2A has the following generalization.

\begin{introprop2B} (Gromov) Equip $S^n$ with a metric $g$, and equip $S^{k_1} \vee S^{k_2}$ with
a metric $h_1$ on $S^{k_1}$ and $h_2$ on $S^{k_2}$.  Suppose that $F: S^n \rightarrow S^{k_1} \vee
S^{k_2}$ has $k_1$-dilation at most $L^{k_1}$.  
Then the linking invariant of $F$ is bounded as follows.

$$|L(F)| \le Iso_{n-k_1}(S^n, g) Vol(g) Area(h_1)^{-1} Area(h_2)^{-1} L^{n+1}.$$

\end{introprop2B}

Similarly, there is a generalization of Proposition 4A using $k_1$-dilation instead
of the Lipschitz constant.  This estimate in turn allows us to bound the k-dilation
of maps from one ellipse to another.

\begin{introprop} Let $E$ be an n-dimensional ellipse
with principal axes $E_0 \le ... \le E_n$, and let $E'$ be an n-dimensional
ellipse with principal axes $E_0' \le ... \le E_n'$.

Suppose that $F$ is a map from $E$ to $E'$
with degree $D$.  Suppose that $k \le
(n+1)/2$.  Let $Q_i$ denote the quotient $E'_i / E_i$.

$$\textrm{Then } [Dil_k(F)]^{\frac{n+1}{k}} \ge c(n) |D| Q_{n+1-k} \prod_{i=1}^n
Q_n .$$

\end{introprop}

This paper contains some new results, but it's also partly expository.
I hope that it will provide some useful background for more 
difficult papers giving estimates for rational homotopy invariants.

One other class of metrics $(S^3, g)$ for which our initial question
has been studied are the Berger spheres.  In \cite{Gr2}, Gromov
gives some estimates for the Hopf invariant of a map from a Berger
sphere.  More generally, he gives estimates for the rational homotopy
invariants of Lipschitz maps from Carnot-Caratheodory spaces.  The
Berger sphere $(S^3, g_t)$ is given by stretching the fibers of the
Hopf fibration to have length $2 \pi t$, while keeping the metric
unchanged in the directions perpendicular to the fibers.  If $t < 1$,
then Proposition 2 gives a good estimate for the Hopf invariant,
but if $t$ is very large then it gives a bad estimate.  In \cite{Gr2},
Gromov proves a better upper bound.  In \cite{G3}, we show that Gromov's
improved upper bound is sharp up to a constant factor.

The
paper \cite{G2} studies second-order rational homotopy
invariants.  The methods in this paper extend immediately to give upper
bounds for second-order homotopy
invariants, but the resulting upper bounds are much too big even for
ellipses.  In \cite{G2}, we give sharp estimates (up to a constant factor)
for the second-order linking invariant of a Lipschitz map from an
n-dimensional ellipse to a wedge of unit spheres.

The paper \cite{G1} studies the Hopf invariant
of continuous maps to a surface of genus $g \ge 2$.  We prove that if $M$ is
an oriented 3-manifold that can be triangulated with $N$ simplices, then
the Hopf invariant of any map from $M$ to a surface of genus 2 is bounded by
$C^N$.  First, we modify
a construction of Milnor-Thurston in order to homotope an arbitrary continuous map
to one with a bound on the 2-dilation.  Then we use the methods from this paper
to bound the Hopf invariant using the 2-dilation.

\vskip5pt

{\bf Acknowledgements.} This paper is based on part of my
thesis.  I thank my thesis advisor, Tom Mrowka, for his help and support.

\section{Gromov's proof}

We begin by proving Gromov's upper bound for the Hopf invariant.

\begin{prop} (Gromov, \cite{Gr1}, pages 96-97) Let $F: (S^3, g)
\rightarrow (S^2, h)$ be a $C^1$ map with 2-dilation $L^2$.  Then the
Hopf invariant of $F$ is at most the following expression.

$$|\textrm{Hopf }(F)| \le Iso_1(g) Vol(g) Area(h)^{-2} L^4.$$

\end{prop}

Remarks: In particular, if the map $F$ has Lipschitz constant $L$,
then it has 2-dilation at most $L^2$ and the conclusion follows.
Also, if $g$ and $h$ happen to be the standard metrics, we get
$|Hopf(F)| \le 100 L^4$.

\proof The first step of the proof is to observe that if $F$ has
small 2-dilation, then one of the fibers $F^{-1}(y)$ has small
length.  This follows from the coarea formula.  We record this
fact as a lemma.  For future reference, we phrase the lemma using
general dimensions.

\begin{smlemma} Suppose that $F: (M^d, g) \rightarrow (N^q, h)$ is
a $C^\infty$ smooth map with q-dilation $\Lambda$, where
$d \ge q$.  Then $F$ has a fiber with (d-q)-dimensional volume at
most $Vol(g) \Lambda / Vol(h)$.
\end{smlemma}

\proof First we write down the coarea formula.

$$\int_M Jac[dF(x)] dvol_g = \int_N Vol^{d-q}[F^{-1}(y)] dvol_h(y).$$

The quantity $Jac[dF(x)]$ is exactly the q-dilation of $F$
at the point $x$, which is at most $\Lambda$.  According to Sard's
theorem, almost every fiber of $F$ is regular.  In particular,
they are almost all manifolds of dimension $d-q$.  The quantity 
$Vol^{d-q}[F^{-1}(y)]$ is the (d-q)-dimensional 
volume of the manifold $F^{-1}(y)$.

It follows that $\int_N Vol^{d-q}[F^{-1}(y)] dvol_h(y) \le \Lambda Vol(M,g)$,
and so we can choose a regular value $y$ so that the fiber $F^{-1}(y)$ has volume
at most $\Lambda Vol(g) / Vol(h)$. \endproof

At this point, we make some remarks about the regularity of $F$.  The coarea
formula holds for any $C^1$ map.  However, for a $C^1$ map, the fibers may
not be generically manifolds.  For a $C^1$ map, the coarea formula involves
the Hausdorff measure of the fibers, and we can find a fiber with (d-q)-dimensional
Hausdorff measure obeying the bound in the lemma.  

On the other hand, any $C^1$ map with q-dilation less than $\Lambda$ may be smoothed
to give a homotopic $C^\infty$ map with q-dilation less than $\Lambda$.  So it
suffices to prove the proposition for $C^\infty$ maps.  From now on we assume
that $F$ is $C^\infty$.

Returning to the proof of Gromov's estimate, we see that $F$ has a
regular fiber $F^{-1}(y)$ with length at most $Vol(g) L^2 / Area(h)$. 
This fiber bounds a 2-chain $C$ with area at most $Vol(g) L^2
Iso_1(g) / Area(h)$.  The map $F$ restricts to a map from $(C,
\partial C)$ to $(S^2, y)$.  This map has a well-defined degree
which is equal to the Hopf invariant of $F$.  On the other hand,
the degree of this map is clearly bounded by $Area(C) L^2 /
Area(h)$.  Plugging in our bound for the area of $C$, we see that
the degree is bounded by $Vol(g) Iso_1(g) L^4 / Area(h)^2$. \endproof

The factor $L^4$ in Gromov's estimate has the right exponent
because of the following example.  (I believe that this example is 
due to Gromov, but it may be older.)

\begin{prop} For large $L$, there is a map $F$ from the unit
3-sphere to the unit 2-sphere with Lipschitz constant $L$ and
Hopf invariant on the order of $L^4$. 
\end{prop}

\proof Let $H: S^3 \rightarrow S^2$ denote the Hopf fibration. 
It has Lipschitz constant $2$.  Let $F: S^2 \rightarrow S^2$
denote a map of Lipschitz constant $L/2$.  For large $L$, the map
$F$ may have degree $D$ on the order of $L^2$.  The composition
$H \circ F$ has Hopf invariant $D^2$ on the order of $L^4$.
\endproof

\section{Linking invariants}

In this section, we recall the definition of the linking invariant.
Gromov's argument directly applies to give estimates for the linking
invariant.

Suppose that $k_1 \le k_2$ and $n = k_1 + k_2 - 1$.  Let $F$ be a
smooth map from $S^n$ to the wedge of spheres $S^{k_1} \vee
S^{k_2}$.  Then we define the linking invariant $L(F)$ to be the
linking number of the fibers $F^{-1}(y_1)$ and $F^{-1}(y_2)$ for
regular points $y_1 \in S^{k_1}$ and $y_2 \in S^{k_2}$.

Here is an equivalent definition. Let $y_1$ be a regular point
in $S^{k_1}$.  Let $C$ be a $k_2$ chain in $S^n$ with boundary
$F^{-1}(y_1)$.  Let
$\pi: S^{k_1} \vee S^{k_2} \rightarrow S^{k_2}$ collapse
$S^{k_1}$ to the basepoint and be equal to the identity on
$S^{k_2}$.  Then $\pi \circ F$ maps $(C, \partial C)$ to
$(S^{k_2}, *)$, where $*$ denotes the basepoint.  We define
$L(F)$ to be the degree of this map.

Standard arguments show that $L(F)$ is independent of the choice
of $y_1$ and $y_2$ and that it is a homotopy invariant.  (These
arguments are completely analogous to the arguments about the Hopf
invariant which appear in the excellent reference \cite{BT}.)

Gromov's argument immediately extends to all linking invariants.

\begin{prop} Suppose $k_1 \le k_2$
and $k_1 + k_2 = n + 1$.  Let $F: (S^n, g)
\rightarrow (S^{k_1} \vee S^{k_2}, h_1 \vee h_2)$
be a map with $k_1$-dilation $L^{k_1}$.  Then the linking
invariant of $F$ is at most the following expression.

$$|L(F)| \le Iso_{k_2 - 1}(g) Vol(g) Area(h_1)^{-1}
Area(h_2)^{-1} L^{n+1}.$$

\end{prop}

Remark on notation: In the statement of the theorem, $h_1$
denotes a metric on $S^{k_1}$ and $h_2$ denotes a metric on
$S^{k_2}$.

\proof Let $U \subset S^n$ denote $F^{-1}(S^{k_1})$.  Clearly the
volume of $U$ is at most the volume of $(S^n, g)$.  Applying
Lemma 1.1 to the map $F: U \rightarrow S^{k_1}$, we find a fiber
$F^{-1}(y_1)$ with $k_2 - 1$-dimensional volume at most $Vol(g)
L^{k_1} / Area(h_1)$. This fiber bounds a $k_2$-dimensional chain
$C$ in $S^n$ of volume at most $Vol(g) L^{k_1} Iso_{k_2 - 1}(g) /
Area(h_1)$.  Now we consider the map $\pi \circ F$ from $(C,
\partial C)$ to $(S^{k_2}, *)$.  The map $\pi$ has Lipschitz
constant 1.  The map $F$ has $k_2$-dilation at most $L^{k_2}$
since $k_2 \ge k_1$.  (This follows from linear algebra, see
Appendix 1 of \cite{G4}.)  Therefore, $\pi \circ F$ has degree at
most $Vol(C) L^{k_2} / Area(h_2)$. Plugging in our bound for
$Vol(C)$ finishes the proof. \endproof

\section{The isoperimetric constants of ellipses}

In this section, we estimate the isoperimetric constants of
ellipses.  This section contains the main work in the
paper.

\begin{prop} If $E$ is an n-dimensional ellipse with principal
axes $E_0 \le ... \le E_n$, then $Iso_k(E) < C(n) 
[E_{k+1} + E_{n-k}]$.
\end{prop}

We will prove in Section 5 that this upper bound is sharp 
up to a constant factor $C(n)$.  The proof uses the k-dimensional
equators of $E$ which we discussed in the introduction.

\proof The ellipse $E$ is $C(n)$-bilipschitz equivalent to the double of a
rectangle $R$ with dimensions $R_1 \le ... \le R_n$, with $R_i =E_i$.  
Each copy of $R$ in the double corresponds to a hemisphere of $E$.

Let $z$ be an integral k-cycle in $E$.  We let $z_N$ denote
the part of $z$ in the Northern hemisphere and $z_S$ denote the
part of $z$ in the Southern hemisphere.  The chain $z_S$ is a
relative cycle in the rectangle $R$.  Our first step is to push
$z_S$ to the boundary of $R$.  More formally, we find a (k+1)-chain
$y_S$ in $R$ with $\partial y_S = z_S + B$ and $B$ in the boundary of
$R$.  The chain $B$ lies
on the equator of $E$, so we can also think of it as a chain in the northern
hemisphere.  Now $z_N - B$ forms a new absolute cycle in the Northern
hemisphere.  The second step is to fill this absolute cycle with a
chain $y_N$.  Finally, we take $y= y_N + y_S$.  It follows immediately
that $\partial y = z_N + z_S = z$.  Our goal is to construct $y_N$ and
$y_S$ with the given bounds on the volume.

The good thing about this approach is that it divides our problem into
two subproblems, and each subproblem takes place in a rectangle
in Euclidean space.  Rectangles are easier to work with than ellipses
in my experience.  This
approach, however, leads immediately to a serious problem: the chain $B$
may be much larger than the cycle $z$.

This point is the main difficulty in the paper.  For example, suppose that $z$ is
the smallest k-dimensional equator.  The intersection of $z$ with each
hemisphere is a k-dimensional rectangle with dimensions $R_1 \times ... \times R_k$
centered in the middle of $R$.  A natural choice for $y$ is given by two copies
of a rectangle with dimensions roughly $R_1 \times ... \times R_{k+1}$.  The situation
is illustrated in Figure 1.

\includegraphics{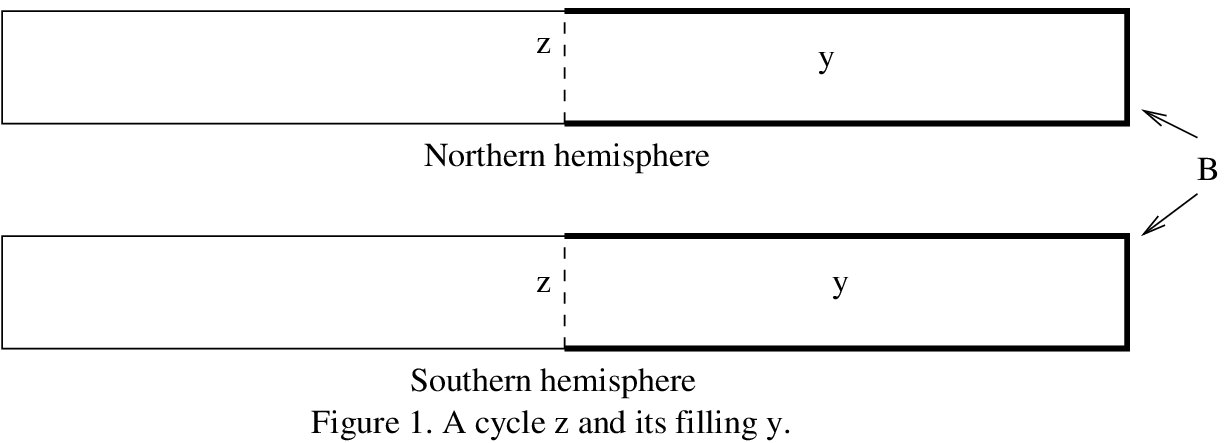}

In the picture, the two dotted lines denote the cycle $z$.  The dotted line in
the Northern hemisphere denotes $z_N$ and the one in the Southern hemisphere
denotes $z_S$.  The chain $y$ is denoted by the labelled rectangles.  
The rectangle in the Northern hemisphere is $y_N$,and the rectangle
in the Southern hemisphere is $y_S$.  Recall that the dimension of $y$ is (k+1) which may
be strictly less than $n$ - so in spite of the picture the reader should not
picture $y$ as an open set.  The picture shows two copies of the chain $B$ -
the chain lies along the equator, and so we can see it in the boundary of each
hemisphere.  The picture shows that $B$ may be much bigger than $z$.

In this example, the cycle $z$ has volume $\sim R_1 ... R_k$, and the chain $B$
has volume $\sim R_2 ... R_{k+1}$.  Consider the last step of our construction:
finding $y_N$.  We need to choose $y_N$ with boundary $z_N - B$.  One
would like to find $y_N$ by applying an isoperimetric inequality for absolute
cycles in $R$.  If the volume of $B$ were comparable to the volume of $z$, then
the isoperimetric inequality would deliver a chain $y_N$ with the desired volume bound.
But as we have seen the volume of $z_N - B$ is much larger than the volume of $z$.  
It is possible to find a cycle $w$ in $R$ with the same volume as $z_N - B$ and
with filling volume much larger than our desired volume bound.  
In order to prove our proposition, we need to take
account of some extra structure in the cycle $z_N - B$, which allows us to fill
it efficiently.  The extra structure comes from keeping track of the "direction"
in which volume is pointing.

(It may perhaps be concluded that trying to build $y$ as a sum of a piece in 
the Northern hemisphere and a piece in the Southern hemisphere is not the right
approach to the isoperimetric problem.  It would be interesting to see another proof.
On the other hand, the directional estimates which we introduce here seem
interesting to me.  In the paper \cite{G2}, they are taken further, leading to
improved results about second-order rational homotopy invariants.)

We begin by defining directional volume.

Let $C$ be an m-chain in $\mathbb{R}^n$.  Suppose that $J$ is an
m-tuple of distinct integers between 1 and $n$.  Let $P(J)$
denote the m-plane with coordinates $x_i$, $i \in J$. We define
the $J$-volume of $C$ to be the volume (with geometric multiplicity) 
of the projection of $C$ to $P(J)$.  For example, if $C$ is the unit
sphere in $\mathbb{R}^n$ and $J$ is any (n-1)-tuple, then the J-volume
of $C$ is twice the volume of the unit (n-1)-ball.  The total volume
of $C$ is roughly equal to the sum of the volumes in different
directions.

$$Vol(C) \le \sum_J Vol_J(C) \le {n \choose m} Vol(C).$$

For an m-tuple $J$, we define $e(J)$ to be one less than the smallest
member of $J$.  For example, if $J$ is $\{ 3, 6, 8 \}$, then $e(J) = 2$.
Since $J$ is an m-tuple of the numbers from $1$ to $n$, $0 \le e(J) \le
n-m$.

In the Northern hemisphere, we will use the following lemma to find
$y_N$.

\begin{lemma} Let $z$ denote an absolute k-cycle in the rectangle $R$.
Then $z$ is the boundary of a chain $y$ obeying the following estimate.

$$|y| \lesssim \sum_{J | e(J) \ge 1} R_{e(J)} Vol_J(z).$$

\end{lemma}

Comments.  If $J$ contains 1, then $e(J) = 0$, and so $Vol_J(z)$ does not appear
on the right-hand side of our inequality.  In other words, the volume of $y$ can
be bounded using only some of the directional volumes of $z$.

\proof Let $\pi(z)$ denote the projection of $z$ to the plane $x_1 = 0$.  The
cylinder from $z$ to $\pi(z)$ has volume at most $R_1 \sum_{e(J) \ge 1} Vol_J(z)$.
The image $\pi(z)$ has the same $J$-volume as $z$ if $J$ does not contain 1, and
it has $J$-volume zero if $J$ contains 1.

Now we repeat this approach, projecting to the cod-2 plane $x_1 = x_2 = 0$, and so
on, until we project to the k-plane $x_1 = ... = x_{n-k} = 0$.  A k-dimensional
cycle inside of this k-plane bounds a (k+1)-chain of volume zero, so after making
these $(n-k)$ projections we are done.

The projection to the plane of codimension $c$, $x_1 = ... = x_c = 0$, costs volume
at most $\sum_{e(J) \ge c} R_c Vol_J(z)$.  Summing the contributions we get the
inequality we want to prove. \endproof

In the Southern hemisphere, we need a slightly more complicated lemma.

\begin{lemma} Let $z$ denote a k-cycle in the rectangle $R$ with $\partial z$
contained in $\partial R$.  Then there is a (k+1)-chain $y$ in $R$ so that
$\partial y = z + B$, where $B$ is contained in $\partial R$, obeying the following
inequalities.

$$|y| \lesssim R_{k+1} |z|.$$

$$\textrm{If $e(J) \ge 1$, } Vol_J(y) = 0.$$

$$\textrm{If $e(J) > 1$, } Vol_J(B) \lesssim |z|. $$

$$\textrm{If $e(J) = 1$, } Vol_J(B) \lesssim (R_{k+1} / R_1 ) |z|. $$

\end{lemma}

\proof The proof is by induction on the dimension $k$ of the cycle.  When $k=0$,
the lemma is trivial.

We consider the intersection of $z$ with planes $x_1 = h$.  We can choose $h$ so
that this intersection is a (k-1)-cycle with volume at most $|z| / R_1$.  We pick
such a value of $h$ and call this intersection $z_h$. 

\includegraphics{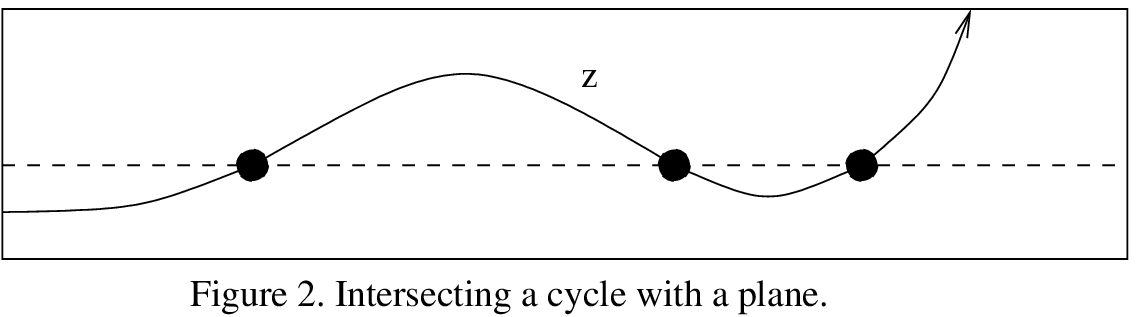}

The solid curve in the figure denotes the relative cycle $z$.  The dotted line
denotes the plane $x_1 = h$.  The three dark points denote their intersection
$z_h$.

We now decompose $z$ into two pieces as follows.

$$z = (z - [0, R_1] \times z_h) + [0, R_1] \times z_h = z_1 + z_2.$$

We deal with the first piece first.  We decompose this first piece as $z_+ + z_-$,
where $z_+$ is the part of $z_1$ lying above the plane $x_1 = h$ and $z_-$ is
the part lying below it.  The two pieces $z_+$ and $z_-$
are each relative cycles.  We illustrate $z_+$ and $z_-$ in the next figure.

\includegraphics{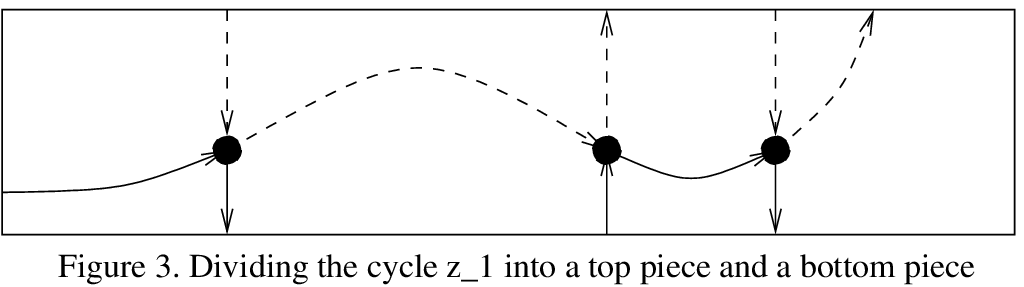}

The dotted curves denote the relative cycle $z_+$ and the solid curves denote
the relative cycle $z_-$.

Now we construct a filling $y_-$ for $z_-$.  In other words, $\partial y_- = z_-
+ B_-$ with $B_- \subset \partial R$.  We construct $y_-$ by pushing $z_-$ down
into the boundary of $R$.  In more detail, we construct a map 
$F: z_- \times [0,1] \rightarrow R$.  The last n-1 coordinates of $F$ are just
the coordinates $x_i$ on $z_+$.  The first coordinate $F_1$ is given by the formula
$F_1(p, t) = t x_1(p)$.  Here $x_1(p)$ denotes the coordinate function of the
point $p \in z_+$.  We define $y_-$ to be the image chain $F([0,1] \times z_-)$.

We illustrate $y_-$ and $B_-$ below.

\includegraphics{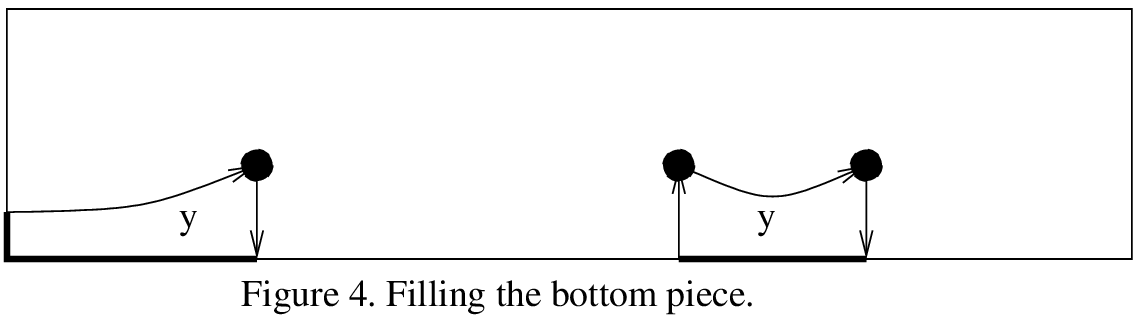}

The chain $y_-$ is labelled $y$ in the figure.  The chain $B$ is the darkened portion
of the boundary of the rectangle.  The solid oriented curves are again $z_-$.

The chain $y_-$ has volume at most $R_1 |z_-| \le 2 R_1 |z|$.  Also, if $J$ is a 
(k+1)-tuple that does not contain $1$, then $Vol_J(y_-) = 0$.  Therefore, $y_-$
obeys the conclusions of the lemma.

We define $B_-$ to be $\partial y_- - z_-$.  Equivalently, $B_-$ is $- F( \{0 \} \times
z_-) + F( [0,1] \times \partial z_-)$.  Since $\partial z_-$ lies in the bottom and 
sides of $R$, not touching the top face of $R$, $B_-$ is contained in $\partial R$.
The first term in $B_-$, $F( \{ 0 \} \times z_-)$, is just the projection of $z_-$
to the plane $x_1 = 0$.  It has volume at most $|z|$.  The second part may have
arbitrarily large volume, but if $J$ does not contain 1, then its $J$-volume is zero.
Therefore $B_-$ obeys the conclusions of the lemma.

Similarly, we define $y_+$ and $B^+$.  Combining them, we get a filling for $z_1$.  Now
we turn to $z_2 = [0,R_1] \times z_h$.  The relative cycle $z_h$ has dimension k-1 
and lives inside an (n-1)-dimensional rectangle with
dimensions $R_2 \times ... \times R_n$.  By induction, we can find a filling $y_h$ for $z_h$
obeying the conclusion of the lemma.

We define $y_2 = [0, R_1] \times y_h$.  The volume of $y_2$ is at most $R_1 |y_h|$.  By
our inductive hypothesis, $|y_h| \le R_{k+1} |z_h|$.  Finally, the volume of $y_2$ is
at most $R_1 R_{k+1} |z_h| \le R_{k+1} |z|$.  Also, if $J$ does not contain 1, then
the $J$-volume of $y_2$ is zero.  So $y_2$ obeys the conclusion of the lemma.

We define $B_2 = \partial y_2 - z_2 = \{ R_1 \} \times y_h - \{0 \} \times y_h + 
[0, R_1] \times B_h.$  By induction, we know that the $J$-volume of $y_h$ vanishes
unless $J$ contains $2$.  Also, if $J$ does not contain 1, then the $J$-volume 
of $[0, R_1] \times B_h$ is zero.  Therefore, if $e(J) > 1$, then the $J$-volume
of $B_2$ is zero.  If $e(J) = 1$, then the $J$-volume of $B_2$ is at most
$2 |y_h| \lesssim R_{k+1} |z_h| \le (R_{k+1}/R_1) |z|$.  Therefore, $B_2$ obeys
the conclusion of the lemma.

Our filling of $z$ is $y= y_1 + y_2$.  We have $\partial y = z + B$ where $B=
B_1 + B_2$.  Because of our estimates for $y_i$ and $B_i$, $y$ and $B$ obey
the conclusions of the lemma. \endproof

Combining the two lemmas, we can finish the proof of the proposition.  Let $z$
be a k-cycle in the ellipse $E$.  Let $z_S$ be the intersection of $z$ with
the Southern hemisphere.  By Lemma 2, we can find a chain $y_S$ in
the Southern hemisphere so that $\partial y_S = z_S + B$, obeying the following
estimates.  First, $|y_S| \lesssim E_{k+1} |z|$.  Second, if $e(J) = 1$, then
$Vol_J(B) \lesssim (E_{k+1} / E_1) |z|$.  Third, if $e(J) > 1$, then $Vol_J(B)
\lesssim |z|$.  Since $B$ is in the equator, we can view it as belonging to
the Northern hemisphere.  Let $z_N$ be the intersection of $z$ with the Northern
hemisphere.  Then $z_N - B$ is an absolute k-cycle in the Northern hemisphere.
We use Lemma 1 to find a chain $y_N$ in the Northern hemisphere with
$\partial y_N = z_N - B$.  According to Lemma 1, the volume of $y_N$
is bounded by the following expression.

$$\sum_{J| e(J) \ge 1} E_{e(J)} Vol_J(z_N - B).$$

If $e(J) > 1$, then $Vol_J(z_N - B) \lesssim |z|$, so the contribution of these
terms is $\lesssim E_{n-k} |z|$.  If $e(J) = 1$, then $Vol_J(z_N - B) \lesssim (E_{k+1} / E_1)
|z|$, and so the contribution of those terms is $\lesssim E_{k+1} |z|$.

Finally, we define $y = y_S + y_N$, a (k+1)-cycle with $\partial y = z$.  The
volume of $y$ is $\lesssim (E_{k+1} + E_{n-k}) |z|$.  Hence $Iso_k(E) \le
C(n) [E_{k+1} + E_{n-k}]$.  \endproof

\section{Upper bounds for homotopy invariants of Lipschitz maps}

In Section 3, we estimated the isoperimetric constants of ellipses.
If we plug these estimates into Propositions 1.1 and 2.1, we get upper bounds
for the Hopf and linking invariants of Lipschitz maps from
ellipses.  

The estimate for the Hopf invariant is as follows.

\begin{prop} If $E$ is a 3-dimensional ellipse and $F: E
\rightarrow (S^2, h)$ is a map with 2-dilation $L^2$ then the
Hopf invariant of $F$ is at most $C(n) E_2 Vol(E) L^4 /
Area(h)^2$.  In particular, if $h$ is the unit sphere metric,
then the Hopf invariant of $F$ is at most $C(n) E_1 E_2^2 E_3 L^4$.
\end{prop}

\proof By Proposition 1.1,

$$|Hopf(F)| \le C Iso_1(E) Vol(E) L^4 Area(h)^{-2}.$$

The volume $Vol(E)$ is at most $C E_1 E_2 E_3$.  By
Proposition 3.1, $Iso_1(E) < C E_2$.  Plugging in we get,

$$|Hopf(F)| \le C L^4 Area(h)^{-2} E_1 E_2^2 E_3.$$

\endproof

The estimate for the linking invariant is as follows.

\begin{prop} Suppose that $E$ is an n-dimensional ellipse with
principal axes $E_0 \le ... \le E_n$.  Suppose that $n +1 = k_1 + k_2$
and $2 \le k_1 \le k_2$.  If $F: E
\rightarrow (S^{k_1} \vee S^{k_2}, h_1 \vee h_2)$ is a map with
$k_1$-dilation $L^{k_1}$, then the linking invariant of $F$ is at
most $C(n) E_{k_2} Vol(E) L^{n+1} Area(h_1)^{-1} Area(h_2)^{-1}$. 
In particular, if each $h_i$ is the unit sphere metric, then the
linking invariant of $F$ is at most $C(n) E_{k_2} Vol(E)
L^{n+1}$.
\end{prop}

\proof By Proposition 2.1, 

$$|L(F)| \le Iso_{k_2 - 1}(E) Vol(E) Area(h_1)^{-1} Area(h_2)^{-1} L^{n+1}.$$

The volume of $E$ is at most $C E_1 ... E_n$.  By Proposition 3.1,
the isoperimetric constant $Iso_{k_2-1}(E)$ is at most $C [E_{k_2} +
E_{n-k_2+1}]$.  This expression is equal to $C [E_{k_2} + E_{k_1}]$.
By assumption, $k_1 \le k_2$, so we conclude the following.

$$|L(F)| \le C L^{n+1} Area(h_1)^{-1} Area(h_2)^{-1} E_{k_2} E_1 ... E_n.$$

\endproof

\section{Lipschitz maps with large Hopf or linking invariants}

The upper bounds in Proposition 4.1 and 4.2 are sharp up to a
constant factor.  To prove this, we construct Lipschitz mappings
with large homotopy invariants.

\begin{prop} If $E$ is a 3-dimensional ellipse and $L > C
E_1^{-1}$, then there exists a map $F$ from $E$ to the unit
2-sphere with Lipschitz constant $L$ and with Hopf invariant at
least $c E_2 Vol(E)$.
\end{prop}

\proof We recall that $E$ is $C$-bilipschitz to the double of the
rectangle $R$ with dimensions $E_1 \times E_2 \times E_3$. 
Inside of the double, we consider a set $U$ equal to the double
of $[0, E_1] \times [E_2 / 3, 2 E_2/3] \times [E_3/ 3, 2 E_3 /
3]$.  Let $\pi$ denote the projection from $U$ to $U' = [E_2 / 3,
2 E_2/3] \times [E_3/ 3, 2 E_3 / 3]$.  Let $f_1$ denote a map
from $(U', \partial U')$ to $(S^2, *)$ with Lipschitz constant
$L$ and degree at least $c E_2 E_3 L^2$.

Let $V$ denote the subset of the Northern hemisphere given by $[0,
E_1] \times \tilde V$, where $\tilde V$ is the region $[0, E_2]
\times [0, E_3] - [E_2/3, 2 E_2 / 3] \times [E_3/3, 2 E_3/3]$. 
The region $V$ is $C$-bilipschitz to a cylinder of the form $S^1(E_3)
\times [0, E_1] \times [0, E_2]$.  Here $S^1(E_3)$ denotes a
circle of radius $E_3$.  Let $\pi$ denote the projection from $V$
to $V' = [0, E_1] \times [0, E_2]$.  
Let $f_2$ denote a map from $(V', \partial V')$ to
$(S^2, *)$ with Lipschitz constant $L$ and degree at least $c E_1
E_2 L^2$.

Finally, we construct a map $F: E \rightarrow S^2$ by combining
$f_1$ and $f_2$.  On the set $U$, we define $F$ to be equal to
the map $f_1 \circ \pi$.  On the set $V$, we define $F$ to be
equal to $f_2 \circ \pi$. On the remainder of $E$, we define $F$
to be the basepoint of $S^2$.

We claim that the map $F$ has Hopf invariant at least $c E_1
E_2^2 E_3 L^4 \sim c E_2 Vol(E) L^4$.  

Topologically, $U, V$ are each thick tubes in $S^3$, and the two
tubes are linked with linking number 1.  For a generic point
$y$ in $S^2$, the fiber
$F^{-1}(y)$ consists of several circles
parallel to the core of $U$ together with several circles
parallel to the core of $V$. More precisely, there are at least
$c E_2 E_3 L^2$ circles parallel to the core of $U$ and at least
$c E_1 E_2 L^2$ circles parallel to the core of $V$.  Now consider
a second generic point $y'$ in $S^2$.  The fiber $F^{-1}(y')$
also consists of $\sim E_2 E_3 L^2$ parallel copies of the core of
$U$ and $\sim E_1 E_2 L^2$ parallel copies of the core of $V$.
The linking number of two parallel copies of the core of $U$ is
zero.  The same holds for $V$.  The linking number of the core of $U$
with the core of $V$ is equal to 1.  Therefore, the linking number
of $F^{-1}(y)$ with $F^{-1}(y')$ is at least $c E_1 E_2^2 E_3 L^4$.
\endproof

\begin{prop} If $E$ is an n-dimensional ellipse and $L > C(n)
E_1^{-1}$, then there exists a map $F$ from $E$ to the wedge of
unit spheres $S^{k_1} \vee S^{k_2}$ with Lipschitz constant $L$
and with linking invariant at least $c(n) E_{k_2} Vol(E)$.
\end{prop}

\proof This proof is essentially the same as the last one, which
boils down to finding two thick linked tubes.

Again we note that $E$ is $C(n)$-bilipschitz to the double of a
rectangle $R$ with dimensions $E_1 \times ... \times E_n$.

We define $U$ to be a thick neighborhood of the double of $[0,
E_1] \times ... \times [0, E_{k_2 - 1}]$.  Topologically, $U$ is
a neighborhood of a standard copy of $S^{k_2 - 1}$ embedded in
$S^n$.  Geometrically, $U$ is $C(n)$-bilipschitz to a product
$U_1 \times [0, E_{k_2}] \times ... \times [0, E_n]$, where $U_1$
is the double of $[0, E_1] \times ... \times [0, E_{k_2 - 1}]$. 
We let $\pi$ denote the projection from $U$ onto $U' = [0,
E_{k_2}] \times ... \times [0, E_n]$.  The dimension of $U'$ is $n-k_2 + 1 =
k_1$.  Let $f_1$ denote a map from $(U', \partial U')$ to
$(S^{k_1}, *)$ with Lipschitz constant $L$ and degree at least
$c(n) E_{k_2} ... E_n L^{k_1}$.

We define $V$ to be a thick neighborhood of the $(k_1 -
1)$-dimensional ellipse $V_1$ sitting in the upper hemisphere,
where $V_1$ is given by the equations $x_1 = ... = x_{k_2 - 1} =
0, \sum_{i=k_2}^n [x_i - (1/2) E_i]^2 E_i^{-2} = 1$.  The region
$V$ is $C(n)$-bilipschitz to $V_1 \times V'$, for $V' = [0, E_1]
\times ... \times [0, E_{k_2}]$.  We let $\pi$ denote the
projection $V \rightarrow V'$.  We let $f_2$ be a map from $(V',
\partial V')$ to $(S^{k_2}, *)$ with Lipschitz constant $L$ and
degree at least $c(n) E_1 ... E_{k_2} L^{k_2}$.

Now we define $F: E \rightarrow S^{k_1} \vee S^{k_2}$.  
We define $F$ on $U$ to be $f_1 \circ \pi$,
which maps $U$ to $S^{k_1} \subset S^{k_1} \vee S^{k_2}$. 
Similarly, we define $F$ on $V$ to be $f_2 \circ \pi$, which maps
$V$ to $S^{k_2} \subset S^{k_1} \vee S^{k_2}$.  On the rest of
$E$, we define $F$ to be the basepoint of $S^{k_1} \vee S^{k_2}$.

The fiber $F^{-1}(y_1)$ for a generic point $y_1 \in S^{k_1}$ is
at least $c(n) E_{k_2} ... E_n L^{k_1}$ parallel copies of the
core of $U$. The fiber $F^{-1}(y_2)$ for a generic point $y_2 \in
S^{k_2}$ is at least $c(n) E_1 ... E_{k_2} L^{k_2}$ parallel
copies of the core of $V$. The core of $U$ and the core of $V$
have linking number 1.  Therefore, the linking invariant of $F$
is at least $c(n) E_{k_2} Vol(E) L^{n+1}$. \endproof

This is also a convenient place to show that our upper bounds for
$Iso_k(E)$ are sharp up to a constant factor.

\begin{prop} The isoperimetric constant $Iso_k(E)$ is at least
$c(n) [E_{k+1} + E_{n-k}]$.
\end{prop}

\proof We let $S \subset E$ denote the smallest
k-dimensional equator.  We let $S'$ denote the largest
(n-k-1)-dimensional equator.  The two spheres $S$ and $S'$
are linked with linking number 1.  We define $V$ to
be a thick neighborhood of $S'$.  The set $V$ is bilipschitz
to $S' \times V'$ for $V' = [0, E_1] \times ... \times [0, E_{k+1}]$.
If $C$ is any chain with $\partial C = S$, then the intersection 
$V \cap C$ has projection that covers $V'$, and so the volume of
$C$ is at least $c(n) E_1 ... E_{k+1}$.  Since the volume of
$S$ is roughly $E_1 ... E_k$, we conclude that $Iso_k(E)$ is at
least $c(n) E_{k+1}$.

Now we let $S \subset E$ denote the smallest (n-k-1)-dimensional
equator.  We let $S'$ denote the largest k-dimensional equator.
The two spheres $S$ and $S'$ have linking number 1.  We define
$U$ to be a thick neighborhood of $S$.  The set $U$ is bilipschitz
to $S \times [0, E_{n-k}] \times ... \times [0, E_n]$.  Therefore,
the filling volume of $S'$ is at least $c(n) E_{n-k} ... E_n$.
Since the volume of $S'$ is roughly $E_{n-k+1} ... E_n$, we conclude
that $Iso_k(E)$ is at least $c(n) E_{n-k}$. \endproof

\section{A bad example}

In this section, we consider a metric more complicated than an
ellipsoidal metric. For this example, the upper bound for
the Hopf invariant in Proposition 1.1 is much bigger than the optimal upper bound.

Our metric $(S^3, g)$ is a connected sum of the unit 3-sphere
with a long thin tube of width $w << 1$ and length
$A >> 1$.  (The long thin tube looks like $S^2(w)
\times [0, A]$, connected on one end to the unit sphere and
capped on the other end with a hemispherical cap of radius $w$.) 
The volume of the tube is roughly $w^2 A$, and we choose $A$ and
$w$ so that $w^2 A$ is much bigger than $1$.  Regardless of $w,
A$, the isoperimetric constant $Iso_1(S^3, g)$ is at least $\sim
1$.  We can consider, for example, a circumference of the unit
3-sphere.  

We consider the Hopf invariants
of Lipschitz maps from $(S^3, g)$ to the unit 2-sphere.  The right-hand
side in Proposition 1.1 is $C Iso_1(g) Vol(g) Area(h)^{-2} L^4$, which
is at least $\sim A w^2 L^4$.  But we will prove that every map
with Lipschitz constant $L$ has Hopf invariant much smaller than
$A w^2 L^4$.

\begin{prop} Every map from $(S^3, g)$ to the unit 2-sphere with 
Lipschitz constant $L$ has Hopf invariant at most $C [1 + A w^3] L^4$.
\end{prop}

\proof Suppose that $F$
is a map from $(S^3, g)$ to the unit 2-sphere with 2-dilation at
most $L^2$.  Define the thick part of $(S^3, g)$ to be the unit
3-sphere together with the first unit length of the connected
tube.  The thin part of $(S^3, g)$ is the rest of it.  Then we
can choose a fiber $F^{-1}(y)$ so that the part of the fiber in
the thick part of $(S^3, g)$ has length at most $L^2$ and the
part of the fiber in the thin part has length at most $L^2 A
w^2$.

Inside the thick part, we included the first unit of the tube,
which is bilipschitz to $S^2(w) \times [0, 1]$.  We can choose a
height $h \in [0, 1]$ so that the fiber $F^{-1}(y)$ meets $S^2(w)
\times \{ h \}$ in at most $L^2$ points.  Then we can perform
surgery at the cut, adding a length at most $L^2 w \le L^2$. 
After the surgery, we have written the fiber as a sum of 1-cycles
$z_1 + z_2$, where $z_1$ lies in the thick part of $(S^3, g)$ and
has length at most $\sim L^2$ and $z_2$ lies in the thin part of
$(S^3, g)$ and has length at most $\sim A w^2 L^2$.

Now we can fill $z_1$ by a 2-chain $C_1$ of area at most $\sim
L^2$.  We can fill $z_2$ more efficiently because it lies in the
thin part of $(S^3, g)$.  Using the argument from Section 3, we
can fill $z_2$ by a 2-chain $C_2$ of area at most $\sim A w^3
L^2$.  Thus we have filled $F^{-1}(y)$ by a chain $C = C_1 + C_2$
of area at most $\sim [1 + A w^3] L^2$.  
Hence the Hopf invariant of $F$ is at most $\sim [1 + A w^3]
L^4$. \endproof

One can choose $A, w$ so that this estimate improves by an
arbitrary factor over the upper bound $A w^2 L^4$.

\section{Application to k-dilation of maps between ellipses}

Our previous results imply a new lower bound for the k-dilation
of a map from one ellipse to another.

\begin{prop} Let $E$ be an n-dimensional ellipse with principal
axes $E_0 \le ... \le E_n$, and let $E'$ be an n-dimensional
ellipse with principal axes $E_0' \le ... \le E_n'$.
Suppose that $F$ is a map from $E$ to $E'$
with degree $D$ and k-dilation $\Lambda$.  Suppose that $k \le
(n+1)/2$.  Let $Q_i$ denote the quotient $E'_i / E_i$.

$$ \Lambda^{\frac{n+1}{k}} \ge c(n) |D| Q_{n+1-k} \prod_{i=1}^n
Q_n .$$

\end{prop}

\proof We use a map from $E'$ with a large linking invariant.  According to
Proposition 5.2, there is a map $\Phi$ from $E'$ to $S^k \vee S^{n+1-k}$ with Lipschitz
constant $L$ and linking invariant at least $c(n) E'_{n+1-k} Vol(E') L^{n+1}$.
Now we consider the composition $\Phi \circ F$ mapping $E$ to $S^k \vee S^{n+1-k}$.
This composition has k-dilation at most $\Lambda L^k$.  Its linking invariant has norm
at least $c(n) |D| E'_{n+1-k} Vol(E') L^{n+1}$.  According to Proposition 4.2, the norm
of the linking invariant must be at most $C(n) E_{n+1-k} Vol(E) L^{n+1} 
\Lambda^{\frac{n+1}{k}}$.  \endproof

\end{document}